\documentclass[11pt]{article}
\usepackage{latexsym,euscript,amsmath,amssymb,amsbsy,amsfonts,amsthm,amsopn,amstext,amsxtra,amscd}
\usepackage{epsfig}
\usepackage{graphics}
\usepackage{url}

\usepackage[english]{babel}

\newtheorem{theorem}{Theorem}

\newcommand{\N}{{\mathbb N}}

\newcommand{\D}{\text{$\mathcal{D}$}}

\begin{document}
\centerline{\bf On statistical independence and}
\centerline{\bf density independence}
\vskip1cm
\centerline{Milan Pa\v st\'eka}
\vskip1cm
{\bf ABSTRACT.}{\it The object of observation in present paper is statistical independence of real sequences and its description as independence with respect to certain class of densities.}
 \footnote{Key words: statistical independence, inependence, asymptotic density, density, uniform distribution, distribution function}
\footnote{Mathematics Subject Classification: 11K06}
\footnote{Research supported by the Grant VEGA 2/0119/23}
\vskip1cm
{\bf Introduction}
The aim of this paper is to characterize statistical independence as independence  of given sequences with respect
a class of selective densities. In first part we recall the notion of statistical independence. The second part is devoted to the definition of selective density and its connection with the distribution of sequences. \newline
{\bf Notation} \newline $\mathcal{X}_S$ - indicator function of the set
$S$ \newline
$v^{-1}(I)= \{n; v(n) \in I\}$ \newline
$|I|$ - length of the interval $I$ \newline
$|A|$ - cardinality of the finite set $A$
\vskip1cm
{\bf 1. Statistical independence.}
The notion of statistically independent sequences is introduced and studied in the work \cite{Ra}.

Let $v_1, \dots, v_m$ be a sequences of elements of a given
interval $[a,b]$, where $-\infty < a < b < \infty$.
The following multi linear forms will play important rule
$$
\Delta_N(h_1, \dots, h_m)=
\frac{1}{N}\sum_{n=1}^N h_1(v_1(n)) \dots h_m(v_m(n)),
$$
$$
\delta_N(h_1, \dots, h_m)=
\prod_{i=1}^m\frac{1}{N}\sum_{n=1}^N h_i(v_i(n)).
$$
for the real functions $h_1, \dots, h_m$ defined on $[a,b]$.
We say that the sequences $v_1, \dots, v_m$ are
{\it statistically independent} if
\begin{equation}
\label{statind}
\lim_{N \to \infty}\Delta_N(f_1, \dots, f_m)- \delta_N(f_1, \dots, f_m)=0,
\end{equation}
for every continuous real functions $f_1, \dots, f_m$ defined
on the interval $[a,b]$. Let us remark that this notion was later studied in various papers. For example \cite{GST}, \cite{GT}.

{\bf 2. Selective density.}
The distribution of sequences was firstly studied by Hermann Weyl in his famous paper \cite{WEY}. We say that a set $S \subset \N$ has asymptotic density if the limit
$$
\lim_{N \to \infty} \frac{|S \cap [1,N]}{N}:=d(S)
$$
exists. In this case the value $d(S)$ is called {\it asymptotic density} of the set $S$. The Weyl's definition of
uniform distribution is equivalent to: given sequence
$\{v(n)\}$ of elements of $[0,1]$ is {\it uniformly distributed mod } $1$ if for each interval $I \subset [0.1]$
the set $v^{-1}(I)$ has asymptotic density and
$d(v^{-1}(I))=1$.

Assume that an increasing sequence of natural numbers.  $\kappa = \{k_1 < k_2 < \dots < k_N < \dots \}$ is given. We say that a set
$S \subset \N$ is $\kappa $ {\it measurable} if the limit
$$
\lim_{N \to \infty} \frac{|S \cap [1, k_N]|}{k_N}:= d_\kappa (S)
$$
exists. In this case the value $d_\kappa(S)$ is called
$\kappa $ density of $S$. The system of all $\kappa$ - measurable sets we shall denote $\D_\kappa$.

A real sequence $v=\{v(n)\}$ will be called $\kappa$ - {\it measurable} if for each real number $x$ the set
$v^{-1}((-\infty, x)$ belongs to $\D_\kappa$. In this case
the function $F(x)= d_\kappa(v^{-1}((-\infty, x))$ we shall
call $\kappa$ - {\it distribution function} of $v$. When
the sequence is bounded and $v(n) \in [a,b]$ for
$-\infty < a < b < \infty$ the the $\kappa$ measurability
of $v$ implies
\begin{equation}
\label{wk}
\lim_{N \to \infty} \frac{1}{k_N}\sum_{n=1}^{k(N)} f(v(n))=
\int_a^b f(x) dF(x)
\end{equation}
for each Riemann - Stjelties integrable function $f$ defined
on the interval $[a,b]$. These types of equalities are known as Weyl citerions. For the proof we refer to \cite{K-N}, \cite{D-T}, \cite{Str}, \cite{Hla}.

Suppose that a sequences $v_i, i=1, \dots, m$ of elements of some interval $[a,b], -\infty < a < b < \infty $ are given. The Helly's selection principe provides the existence of such
sequence of natural numbers $\kappa$ that these sequences are
$\kappa$ - measurable. We say that these sequences are $\kappa$ - independent if the set
$v_1^{-1}((a,x_1))\cap \dots \cap v_m^{-1}((a,x_m))$ belongs to
$\D_\kappa$ and
\begin{equation}
\label{kappai}
d_\kappa(v_1^{-1}([a,x_1))\cap \dots \cap v_m^{-1}([a,x_m)))=
F_1(x_1) \dots F_m(x_m),
\end{equation}
for every $x_1, \dots, x_m \in [a,b]$ such that $F_i$ is continuous in the point $x_i, i=1, \dots, m$. In this circumstances is interesting the following equality.
Put
$$
v_1^{-1}([a,x_1))\cap \dots \cap v_m^{-1}([a,x_m))=S(x_1, \dots, x_m):=S.
$$
Taking account that $\mathcal{X}_I(v(n))= \mathcal{X}_{v^{-1}(I)}(n)$ for arbitrary interval $I$ and
real sequence $\{v(n)\}$ we get
\begin{equation}
\label{intersection}
\Delta_N(\mathcal{X}_{[a,x_1)}, \dots, \mathcal{X}_{[a,x_m)})=
\frac{1}{N}\sum_{n=1}^N \mathcal{X}_S.
\end{equation}
Let us remark that for $\kappa$ measurable
If the sequences $v_i, i=1, \dots, m$ we have
\begin{equation}
\label{produktintegral}
\lim_{N \to \infty} \delta_{k_N}(h_1, \dots, h_m)=
\prod_{i=1}^m \int_a^b h_i(x)dF_i(x),
\end{equation}
where $F_i$ is $d_\kappa$ - distribution function of
$v_i$ and $h_i$ are Riemann Stjelties integrable
functions with respect to $F_i$.
\vskip1cm
{\bf 3. The equivalence.} Now we can formulate the following assertion:
\begin{theorem} Let $\{v_1(n)\}, \{v_2(n)\}, \dots, \{v_m(n)\}$ be bounded sequences of real numbers. Then the following statements are eqiuvalent: \newline
(i) given sequences are statistically independent \newline
(ii) given sequences are  $\kappa$ - independent for each sequence $\kappa$ such that these sequence are $\kappa$ measurable.
\end{theorem}

{\bf Proof.} (i)$\Rightarrow$(ii). Suppose that given sequences are statistically independent.  Let $x_i$ is a point of continuity of $F_i$, $i=1, \dots, m$. Then for given $\varepsilon>0$
 such continuous functions $0 \le f_i \le \overline{f}_i$, $i=1, \dots, m$ exist that
\begin{equation}
\label{chinerovnost}
f_i < \mathcal{X}_{[a,x_i)} < \overline{f}_i
\end{equation}
and
\begin{equation}
\label{chiepsilon}
\int_a^b\overline{f}_i(x)dF_i(x) - \int_a^bf_i(x)dF_i(x) \le \varepsilon.
\end{equation}
Taking account
$\mathcal{X}_{[a,x_i)}(v_i(n))=\mathcal{X}_{v_i^{-1}[a,x_i)}(n)$
we get from (\ref{chinerovnost})
$$
f_i(v_i(n)) < \mathcal{X}_{[a,x_i)}(n) < \overline{f}_i(v_i(n))
$$
for $n \in \N$ and $i=1, \dots, m$. This yields
$$
\Delta_N(f_1, \dots, f_m)\le
\Delta_N(\mathcal{X}_{[a,x_1)}, \dots, \mathcal{X}_{[a,x_m)})\le \Delta_N(\overline{f}_1, \dots, \overline{f}_m)
$$
and so (\ref{intersection}) yields
\begin{equation}
\label{Delta}
\Delta_N(f_1, \dots, f_m)\le
\frac{1}{N}\sum_{n=1}^N \mathcal{X}_S\le \Delta_N(\overline{f}_1, \dots, \overline{f}_m).
\end{equation}
Moreover it holds
\begin{equation}
\label{delta}
\delta_N(f_1, \dots, f_m)\le
\delta_N(\mathcal{X}_{[a,x_1)}, \dots, \mathcal{X}_{[a,x_m)})\le \delta_N(\overline{f}_1, \dots, \overline{f}_m).
\end{equation}
Suppose that $\kappa =\{k(N)\}$ is such increasing sequence
of natural numbers that the sequences $v_i, i=1, \dots, m$ are
$\kappa$ measurable. Then
$$
\lim_{N \to \infty} \delta_{k(N)}
(\mathcal{X}_{[a,x_1)}, \dots, \mathcal{X}_{[a,x_m)})=
\prod_{i=1}^m d_\kappa(v_i^{-1}([a,x_i))).
$$
We suppose that the sequences $\{v_1(n)\}, \dots, \{v_m(n)\}$ are statistically independent and so we can apply (\ref{statind}) for $k(N)$ and we get

$$
\lim_{N \to \infty} \Delta_{k(N)}(f_1, \dots, f_m)=
\prod_{i=1}^m \int_a^b f_i(x)dF_i(x)$$ and \newline
$$\lim_{N \to \infty} \Delta_{k(N)}(\overline{f}_1, \dots, \overline{f}_m)=
\prod_{i=1}^m \int_a^b \overline{f}_i(x)dF_i(x)$$.
And so
$$
\prod_{i=1}^m \int_a^b f_i(x)dF_i(x) \le \underline{\lim_{N \to \infty}}\frac{1}{k(N)}\sum_{n=1}^{k(N)} \mathcal{X}_S\le \prod_{i=1}^m \int_a^b \overline{f}_i(x)dF_i(x)
$$
$$
\prod_{i=1}^m \int_a^b f_i(x)dF_i(x) \le \overline{\lim_{N \to \infty}}\frac{1}{k(N)}\sum_{n=1}^{k(N)} \mathcal{X}_S\le \prod_{i=1}^m \int_a^b \overline{f}_i(x)dF_i(x).
$$
Clearly
$$
\prod_{i=1}^m \int_a^b f_i(x)dF_i(x) \le \prod_{i=1}^m
d_{\kappa}(v_i^{-1}([a, x_i)) \le \prod_{i=1}^m \int_a^b \overline{f}_i(x)dF_i(x).
$$
And so taking account (\ref{chiepsilon}) and
$\varepsilon \to 0^+$ we get
$$
\prod_{i=1}^m
d_{\kappa}(v_i^{-1}([a, x_i))=
\lim_{N \to \infty} \frac{1}{k(N)}\sum_{n =1}^{k(N)}
\mathcal{X}_S = d_{\kappa}(S).
$$

 (ii)$\Rightarrow$ (i). Let $\mathcal{S}$ be the set of all step functions
$s= a_1\mathcal{X}_{I_1} +\dots + a_\ell\mathcal{X}_{I_\ell}$ where the endpoints of the intervals $I_j$ are points of continuity of $F_1, \dots, F_m$. Then for every $g_1, \dots, g_m \in \mathcal{S}$ we have
$$
\lim_{N \to \infty} \Delta_{k(N)}(g_1, \dots, g_m) - \delta_{k(N)}(g_1, \dots, g_m) =0.
$$
And so from (\ref{produktintegral}) we get
\begin{equation}
\label{z1}
\lim_{N \to \infty} \Delta_{k(N)}(g_1, \dots, g_m)= \int_a^bg_1(x)dF_1(x)\dots \int_a^bg_m(x)dF_m(x)
\end{equation}
Let us consider a real functions $f_1, \dots, f_m$ continuous on the interval $[a,b]$. These functions are bounded, thus we can suppose to be positive. For each $\varepsilon >0$  such step functions
$s_1, \dots, s_m, S_1, \dots, S_m \in \mathcal{S}$ exist that
\begin{equation}
\label{nerov}
s_i \le f_i \le S_i,
\end{equation}
\begin{equation}
\label{epsilon}
\int_a^b S_i(x) dF_i(x) -
\int_a^b s_i(x) dF_i(x) < \varepsilon,
\end{equation}
for $i=1, \dots, m$. This yields
$$
\int^a_b s_i(x)dF_i(x) \le \int_a^bf_i(x)dF_i(x) \le \int_a^b S_i(x)dF_i(x).
$$
Denote for simplicity
$$
\Delta_N=\Delta_N(f_1, \dots, f_m), \
\delta_N= \delta_N(f_1, \dots, f_m).
$$
Then we have
$$
\Delta_N(s_1, \dots, s_m)\le \Delta_N \le
\Delta_N(S_1, \dots, S_m),
$$
and
$$
\delta_N(s_1, \dots, s_m)\le \delta_N \le
\delta_N(S_1, \dots, S_m).
$$
Applying (\ref{z1}) to the step functions $s_1, \dots, s_m$ and  $S_1, \dots, S_m$ we get
$$
\prod_{i=1}^m \int_a^b s_i(x)dF_i(x)
\le \underline{\lim}_{N\to \infty} \Delta_{k(N)}\le
\overline{\lim}_{N\to \infty} \Delta_{k(N)}\le \prod_{i=1}^m \int_a^b S_i(x)dF_i(x).
$$
Clearly
$$
\prod_{i=1}^m \int_a^b s_i(x)dF_i(x)
\le \lim_{N \to \infty} \delta_{k(N)}\le \prod_{i=1}^m \int_a^b S_i(x)dF_i(x).
$$
This yields
$$
|\underline{\lim_{N\to \infty}} \Delta_{k(N)}- \lim_{N \to \infty} \delta_{k(N)}|  \le
\prod_{i=1}^m \int_a^b S_i(x)dF_i(x)-\prod_{i=1}^m \int_a^b s_i(x)dF_i(x).
$$
Analogously
$$
|\overline{\lim_{N\to \infty}} \Delta_{k(N)}- \lim_{N \to \infty} \delta_{k(N)}|  \le
\prod_{i=1}^m \int_a^b S_i(x)dF_i(x)-\prod_{i=1}^m \int_a^b s_i(x)dF_i(x).
$$
Taking account the  inequalities (\ref{epsilon})
and $\varepsilon\to 0^+$ we get
\begin{equation}
\label{LS}
\lim_{N \to \infty}\Delta_{k(N)}-\delta_{k(N)}=0.
\end{equation}
Each sequence of natural numbers contains by Helly's selection principe such subsequence $\kappa$ that the sequences $\{v_1(n)\}, \dots, \{v_m(n)\}$ are $\kappa$ measurable. Therefore from (\ref{LS}) we get, that
$\overline{\lim}_{N \to \infty}|\Delta_{N}-\delta_{N}|=0.$
\qed \newline
Let us remark that this proof is similar to the proof of
Theorem 3 in \cite{P-T1}.

Author's adress: Department of Mathematics and Informatics, Faculty of Education, University of Trnava,
Priemyselna 4, Trnava, Slovakia.

\end{document}